\documentclass[12pt]{article}
\usepackage{amssymb}
\usepackage{latexsym,bm}
\usepackage{graphicx}
\usepackage{amsmath}
\usepackage{mathrsfs}
\usepackage{epstopdf}

\date{}
\newcounter{mathitem}
  {\begin{list}{{$(\roman{mathitem})$}}{
   \setcounter{mathitem}{0}
   \usecounter{mathitem}
   \setlength{\topsep}{0pt plus 2pt minus 0pt}
   \setlength{\parskip}{0pt plus 2pt minus 0pt}
   \setlength{\partopsep}{0pt plus 2pt minus 0pt}
   \setlength{\parsep}{0pt plus 2pt minus 0pt}
   \setlength{\leftmargin}{35pt}
   \setlength{\itemsep}{0pt plus 2pt minus 0pt}}}
  {\end{list}}

\begin{document}
\title{The detour covering number and cummerbund covering number of a graph
\footnote{E-mail addresses: {\tt lichengli0130@126.com} (C. Li), {\tt zhan@math.ecnu.edu.cn} (X. Zhan).}}
\author{\hskip -10mm Chengli Li and Xingzhi Zhan\thanks{Corresponding author}\\
{\hskip -10mm \small Department of Mathematics,  Key Laboratory of MEA (Ministry of Education) }\\
{\hskip -10mm \small \& Shanghai Key Laboratory of PMMP, East China Normal University, Shanghai 200241, China}}\maketitle
\begin{abstract}
We introduce several new concepts about graphs and investigate their basic properties. A longest path in a graph is called a detour and a longest cycle
is called a cummerbund. The detour covering number of a graph is the number of vertices that lie in a detour. A graph is said to be detour covered
if every vertex lies in a detour. The cummerbund covering number and cummerbund covered graphs are defined similarly. Some of the main results are as follows.
(1) Minimum degree and forbidden subgraph conditions that ensure a graph to be cummerbund covered or detour covered. (2) The minimum cummerbund covering number
and minimum detour covering number of a graph with connectivity or girth conditions. (3) The minimum cummerbund covering number of a $2$-connected
bipartite graph and the extremal graphs.
\end{abstract}


\vskip 8mm

\section{Introduction}

We consider finite simple graphs and use standard terminology and notation from [4] and [13]. The {\it order} of a graph is its number of vertices, and the
{\it size} its number of edges. In this paper we will introduce several new concepts about graphs and investigate their basic properties.

{\bf Definition 1.} A longest path in a graph is called a {\it detour} and a longest cycle is called a {\it cummerbund}.

It seems that the concise term detour first appeared in [11] and now it has been widely used (e.g. [2] and [5]). Professor Roger Horn of the University of Utah
suggested that a longest cycle be called a cummerbund since a shortest cycle is called a girdle [3, p.117].

{\bf Definition 2.} The {\it detour covering number} of a graph $G,$ denoted by ${\rm dc}(G),$ is the number of vertices of $G$ that lie in a detour, and the {\it cummerbund covering number} of a graph $G,$ denoted by ${\rm cc}(G),$ is the number of vertices of $G$ that lie in a cummerbund.

{\bf Definition 3.} A graph $G$ is said to be {\it detour covered} if every vertex of $G$ lies in a detour, and a graph $G$ is said to be {\it cummerbund covered}
if every vertex of $G$ lies in a cummerbund.

Clearly, cummerbund covered graphs are a generalization of hamiltonian graphs. These new concepts might be useful in modeling real-world problems. Imagine a city
as a graph $G$ in which vertices represent the people and cummerbunds represent the clubs. A person represented by the vertex $v$ belongs to a club represented by the cummerbund
$C$ if and only if $v$ lies in $C.$ Then the cummerbund covering number of $G$ means the number of those people who belong to some club, and the statement that everybody in the city belongs to some club is described as $G$ being cummerbund covered.

In Section 2 we state and prove the main results.  Some of the main results are as follows.
(1) Minimum degree and forbidden subgraph conditions that ensure a graph to be cummerbund covered or detour covered. (2) The minimum cummerbund covering number
and minimum detour covering number of a graph with connectivity or girth conditions. (3) The minimum cummerbund covering number of a $2$-connected
bipartite graph and the extremal graphs. Finally we pose several unsolved problems.

We denote by $V(G)$ and $E(G)$ the vertex set and edge set of a graph $G,$ respectively, and denote by $|G|$ and $||G||$ the order and size of $G,$ respectively. Thus, if $H$ is a cycle or a path then $||H||$ means the length of $H.$ The neighborhood and degree of a vertex $x$ in a graph $G$ is denoted by $N_G(x)$ and $\deg_G(x),$ respectively.
If the graph $G$ is clear from the context, the subscript $G$ might be omitted.  We denote by $\delta(G)$ the minimum degree of $G.$ For a vertex subset $S\subseteq V(G),$ we use $G[S]$ to denote the subgraph of $G$ induced by $S.$ Given two disjoint vertex subsets $S$ and $T$ of $G,$ we denote by $[S,T]$ the set of edges having one endpoint in $S$ and the other in $T.$ An $(S,T)$-path is a path which starts at a vertex of $S,$ ends at a vertex of $T,$ and whose internal vertices belong to neither $S$ nor $T.$ In particular, if $S=\{x\}$, we write an $(\{x\},T)$-path as an $(x,T)$-path. An $(x,y)$-path is a path with endpoints $x$ and $y.$

We denote by $C_n$ and $K_n$ the cycle of order $n$ and the complete graph of order $n,$ respectively. As usual, $qK_2$ denotes the graph consisting of $q$ pairwise vertex-disjoint edges. For two graphs $G$ and $H,$ $G\vee H$ denotes the join of $G$ and $H,$ which is obtained from the disjoint union $G+H$ by adding edges joining every vertex of $G$ to every vertex of $H.$ The distance of two vertices $u$ and $v$ in a graph is denoted by $d(u,v).$ We denote by $c(G),$ $g(G),$ and $\kappa(G)$ the circumference, girth, and connectivity of $G$, respectively.

For graphs we will use equality up to isomorphism, so $G=H$ means that $G$ and $H$ are isomorphic.

\section{Main results}

The {\it detour order} of a graph is the order of a detour in the graph; i.e., the number of vertices in a detour.
We first observe that neither of the class of detour covered graphs and the class of cummerbund covered graphs is contained in the other.

{\bf Definition 4.} {\it Duplicating a vertex $v$} in a graph is the operation of adding a new vertex $v^{\prime}$ and adding edges incident to $v^{\prime}$  such that  the neighborhood of $v^{\prime}$
is  $N(v).$

{\bf Theorem 1.} {\it For every integer $n\ge 7,$ there exists a $2$-connected detour covered graph of order $n$ which is not cummerbund covered. For every integer $n\ge 12,$
there exists a $2$-connected cummerbund covered graph of order $n$ which is not detour covered.}

{\bf Proof.} For every $n\ge 7,$ we construct a $2$-connected detour covered graph $G_n$ of order $n$ which is not cummerbund covered. First let $G_7$ be the graph
in Figure 1(a).
\begin{figure}[h]
\centering
\includegraphics[width=0.7\textwidth]{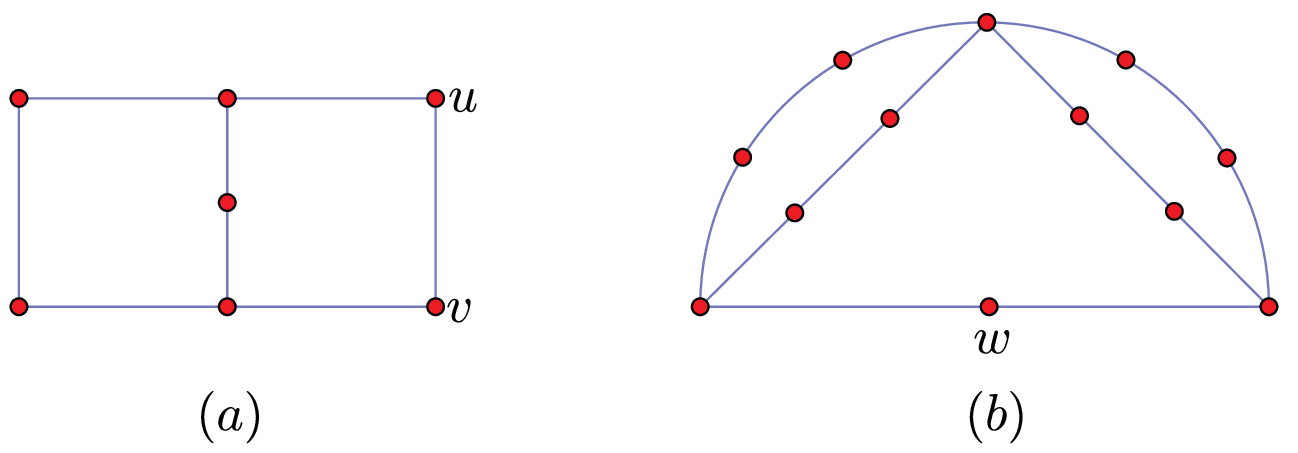}
\caption{The graphs $G_7$ and $H_{12}$}
\end{figure}

Then for $n\ge 8,$ $G_n$ is obtained from $G_7$ by successively subdividing the edge $uv$ $n-7$ times. The detour order of $G_n$ is $n;$ i.e., it is traceable.
The circumference of $G_n$ is $n-1.$

For every $n\ge 12,$ we construct a $2$-connected cummerbund covered graph $H_n$ which is not detour covered. $H_{12}$ is depicted in Figure 1(b). For $n\ge 13,$
$H_n$ is obtained from $H_{12}$ by duplicating the vertex $w$ $n-12$ times. Clearly, the circumference and detour order of $H_n$ are $8$ and $11,$ respectively. $H_n$ has
four detours and $4(n-11)$ cummerbunds. It is easy to verify that $G_n$ and $H_n$ satisfy the requirements. \hfill $\Box$

An {\it empty graph} is a graph without edges, while the {\it null graph} is the graph without vertices. Thus, the null graph is an empty graph. Let $H$ be a subgraph of
a graph $G.$ $H$ is said to be {\it dominating} in $G$ if $G-V(H)$ is an empty graph. Very often, $H$ is a cycle or a path.

Dirac's classic theorem [6] (see also [4, p.485] or [13, p.288]) states that every graph of order $n\ge 3$ with minimum degree at least $n/2$ is hamiltonian. Next in Theorem 3
we prove its counterpart for cummerbund covered graphs.  We will need the following elegant result of Hoa [10].

{\bf Lemma 2 [10].} {\it Let $G$ be a $2$-connected graph of order $n$ with $\delta(G)\ge n/3.$ If $C$ is a cummerbund of $G,$ then $G-V(C)$ is either an empty graph or
a complete graph.}

{\bf Theorem 3.} {\it If $G$ is a $2$-connected graph of order $n$ with $\delta(G)\ge n/3,$ then $G$ is cummerbund covered.}

{\bf Proof.} To the contrary, suppose that $G$ is a $2$-connected graph of order $n$ with $\delta(G)\ge n/3,$ but  $G$ is not cummerbund covered.
Then $G$ contains a vertex $x$ which does not lie in any cummerbund. Let $C$ be a cummerbund of $G.$ We have $x\not\in V(C).$ By Lemma 2, $G-V(C)$ is either an empty graph or a complete graph. We distinguish two cases.

Case 1.  $G-V(C)$ is an empty graph.

In this case $N_G(x)=N_C(x),$ and $|N_C(x)|\ge 2$ since $G$ is $2$-connected. We assert that for any distinct $u, v\in N_C(x),$ $d_C(u,v)\ge 3.$
In fact, if $d_C(u,v)=1,$ then $G$ contains a cycle longer than $C,$ a contradiction, and if $d_C(u,v)=2,$ then $x$ lies in a cummerbund, a contradiction again.
The condition $\delta(G)\ge n/3$ implies that $|N_C(x)|\ge n/3.$ We deduce that
$$
n\ge 1+||C||\ge 1+3\cdot \frac{n}{3}=n+1,
$$
a contradiction.

Case 2.  $G-V(C)$ is a complete graph.

Denote $H=G-V(C)$ and let $|H|=p.$ Note that if $p=1,$ then $H$ is an empty graph, which has been treated in Case 1. Next assume $p\ge 2.$
Then $H=K_p$ and $|C|=n-p.$ Note that $x\in V(H).$ Choose a vertex $y\in V(H)$ such that $y\neq x.$ Let $N_C(x)\cup N_C(y)=\{w_1,\ldots,w_s\}$
 where $w_1,\ldots,w_s$ appear on $C$ in order and we set $w_{s+1}=w_1.$  The vertices $w_1,\ldots,w_s$ divide the cycle $C$ into $s$ paths $C(w_i, w_{i+1}),$ $1\le i\le s.$ A path $C(w_i, w_{i+1})$ is called an {\it exceptional path of $C$} if its length $||C(w_i, w_{i+1})||\ge p+2.$ Let $k$ be the number of exceptional paths of $C.$ Dirac's theorem
 [6, Theorem 4] states that the circumference of a $2$-connected graph of order $n$ and minimum degree $d$ is at least ${\rm min}\{2d, n\}.$ Thus $|C|\ge 2n/3$ and hence
 $|H|=p\le n/3.$ It follows that every vertex of $H$ has at least one neighbor on $C,$ since $\delta(G)\ge n/3.$

Subcase 2.1. $N_C(x)\neq N_C(y).$

We will repeatedly use the fact that if $w_i\in N_C(x)$ and $w_{i+1}\in N_C(y)$ or $w_{i+1}\in N_C(x)$ and $w_i\in N_C(y)$ then $C(w_i, w_{i+1})$ is an exceptional path of $C,$ since $C$ is a cummerbund and $x$ does not lie in any cummerbund. It is easy to verify that
$s\ge k\ge {\rm max}\{2,\, 1+|N_C(x)\cap N_C(y)|\}.$ Since $C$ is a cummerbund, we have $d_C(w_i, w_{i+1})\ge 2$ for each $i.$ Hence
$$
||C||\ge k(p+2)+2(s-k)=kp+2s. \eqno (1)
$$
Note that
$$
d_G(x)-d_C(x)=p-1=d_G(y)-d_C(y), \eqno (2)
$$
and $s=|N_C(x)\cup N_C(y)|=d_C(x)+d_C(y)-|N_C(x)\cap N_C(y)|.$ Using (2) we obtain
$$
s=d_G(x)+d_G(y)-2(p-1)-|N_C(x)\cap N_C(y)|.  \eqno (3)
$$
Since $k\ge {\rm max}\{2,\, 1+|N_C(x)\cap N_C(y)|\}$ and $p\ge 2,$ we have $(k-2)p+2\ge 2k-2\ge k\ge 1+|N_C(x)\cap N_C(y)|;$ i.e.,
$$
(k-2)p+1-|N_C(x)\cap N_C(y)|\ge 0. \eqno (4)
$$
Using (1), (2), (3), (4) and the fact that $s\ge d_C(x)$ we obtain
\begin{align*}
n&\ge ||C||+p\ge (kp+s)+s+p\\
 &=(kp+s)+[d_G(x)+d_G(y)-2(p-1)-|N_C(x)\cap N_C(y)|]+[d_G(x)-d_C(x)+1]\\
 &=2d_G(x)+d_G(y)+[s-d_C(x)]+2+[(k-2)p+1-|N_C(x)\cap N_C(y)|]\\
 &\ge 2d_G(x)+d_G(y)+2\\
 &\ge 2\cdot \frac{n}{3}+\frac{n}{3}+2\\
 &\ge n+2,
\end{align*}
which is a contradiction.

Subcase 2.2. $N_C(x)=N_C(y).$

Now $s=k=d_C(x)=d_G(x)-p+1\ge n/3-p+1.$ Thus
$$
n=||C||+p\ge s(p+2)+p\ge (n/3-p+1)(p+2)+p\triangleq f(p). \eqno (5)
$$
Since $C$ is a cummerbund, $x$ does not lie in any cummerbund and $G$ is $2$-connected, considering a cycle containing all the vertices of $H$ and a path on $C$ we deduce that
$2[(p-1)+2]\le n-p-1;$ i.e., $2\le p\le (n-3)/3.$ Observe that the minimum value of $f(p)$ can be attained at either $2$ or $(n-3)/3.$ By (5) and the fact that $n\ge 9$ we deduce that
$n\ge {\rm min}\{f(2), f((n-3)/3)\}=n+1,$ a contradiction. This completes the proof. \hfill $\Box$

{\bf Remark 1.} The lower bound $n/3$ for the minimum degree in Theorem 3 is best possible. Since $\delta(G)$ is an integer, the condition $\delta(G)\ge n/3$ is equivalent to
$\delta(G)\ge \lceil n/3\rceil.$ For every integer $n\ge 9,$ we construct a $2$-connected graph $H_n$ of order $n$ with $\delta(H_n)=\lceil n/3\rceil -1$ which is not cummerbund
covered. If $n=3k,$ let $H_n=(k-1)K_1\vee ((k-1)K_2+3K_1);$ if $n=3k+1,$ let $H_n=kK_1\vee (kK_2+K_1)$ and if $n=3k+2,$ let $H_n=kK_1\vee (kK_2+2K_1).$

{\bf Remark 2.} Using similar arguments, Theorem 3 can also be proved by applying a result of  Bauer, Schmeichel and Veldman [1, Theorem 7].

Recall that ${\rm dc}(G)$ and ${\rm cc}(G)$ denote the detour covering number and cummerbund covering number of a graph $G,$ respectively.
The following lemma is easy to verify.

{\bf Lemma 4.} {\it Let $G$ be a graph and let $H=G\vee K_1$ where $V(K_1)=\{x\}.$ Then $\kappa(G)=k$ if and only if $\kappa(H)=k+1.$ Suppose $P$ is a $(y,z)$-path in $G.$
Then $P$ is a detour of $G$ if and only if $P\cup\{xy, xz\}$ is a cummerbund of $H.$ Consequently, ${\rm dc}(G)={\rm cc}(H)-1.$ In particular, $G$ is detour covered if and
only if $H$ is cummerbund covered.}

{\bf Theorem 5.} {\it If $G$ is a connected graph of order $n$ with $\delta(G)\ge (n-2)/3,$ then $G$ is detour covered.}

{\bf Proof.} The conclusion follows from Lemma 4 and Theorem 3. \hfill $\Box$

{\bf Remark 3.} The lower bound $(n-2)/3$ for the minimum degree in Theorem 5 is best possible. Since $\delta(G)$ is an integer, the condition $\delta(G)\ge (n-2)/3$ is
equivalent to $\delta(G)\ge \lceil (n-2)/3\rceil.$ For every integer $n\ge 6,$ we construct a connected graph $R_n$ of order $n$ with $\delta(R_n)=\lceil (n-2)/3\rceil -1$ which is not detour covered. If $n=3k,$ let $R_n=(k-1)K_1\vee (kK_2+K_1);$ if $n=3k+1,$ let $R_n=(k-1)K_1\vee (kK_2+2K_1)$ and if $n=3k+2,$ let $R_n=(k-1)K_1\vee (kK_2+3K_1).$

Recall that $c(G)$ denotes the circumference of a graph $G.$

We will need the following result of Fang and Xiong [8, the case $m=0$ of Theorem 8].

{\bf Lemma 6 [8].} {\it Let $G$ be a graph with $\kappa(G)=k\ge 2.$ If $c(G)\le 3k-1,$ then every cummerbund of $G$ is dominating.
}

{\bf Theorem 7.} {\it The minimum cummerbund covering number of a $k$-connected graph of order $n$ where $k\ge 2$ is $\min \{n,3k\}.$ }

{\bf Proof.} Let $G$ be a $k$-connected graph of order $n.$ Then $\delta(G)\ge k.$ We first show that ${\rm cc}(G)\ge\min \{n,3k\}.$

 If $n\le 3k,$ then Theorem 3 implies that ${\rm cc}(G)=n=\min \{n,3k\}.$

Next we assume that $n>3k.$ If $c(G)\ge 3k,$ then ${\rm cc}(G)\ge 3k=\min \{n,3k\}.$ Suppose that $c(G)\le 3k-1.$ Now $G$ is nonhamiltonian.
Let $C$ be a cummerbund of $G.$ By Lemma 6, $C$ is a dominating cycle of $G.$ Let $v$ be an arbitrary but fixed vertex in $V(G)\setminus V(C).$
Since $c(G)\le 3k-1$ and $\deg_C(v)= \deg_G(v)\ge k,$ by the pigeonhole principle, there exist two vertices $x,y\in V(C)$ such that $vx,vy\in E(G)$ and $d_C(x,y)\le 2.$
If $d_C(x,y)=1,$ then there exists a cycle longer than $C$, a contradiction. Hence, $d_C(x,y)=2,$ and so $v$ lies on a cummerbund. Since $v$ was arbitrarily chosen,
$G$ is cummerbund covered. It follows that ${\rm cc}(G)=n> 3k=\min \{n,3k\}.$

Conversely we give examples to show that the lower bound $\min \{n,3k\}$ can be attained. If $n\le 3k,$ by Theorem 3, any $k$-connected graph of order $n$ attains $n=\min \{n,3k\}.$
If $n>3k,$ the graph $kK_1\vee (kK_2+(n-3k)K_1)$ is a $k$-connected graph of order $n$ and circumference $3k$ which has cummerbund covering number $3k=\min \{n,3k\}.$
This completes the proof of Theorem 7. \hfill$\Box$

{\bf Theorem 8.} {\it The minimum detour covering number of a $k$-connected graph of order $n$ is $\min \{n,3k+2\}.$ }

{\bf Proof.} Let $G$ be a $k$-connected graph of order $n.$ Construct a graph $H=G\vee K_1.$ Then $H$ is a $(k+1)$-connected graph of order $n+1.$
By Lemma 4 and Theorem 7,
$$
{\rm dc}(G)={\rm cc}(H)-1\ge \min \{n+1,3(k+1)\}-1=\min \{n,3k+2\}.
$$

Conversely, we give examples to show that the lower bound $\min \{n,3k+2\}$ can be attained. If $n\le 3k+2,$ by Theorem 5, any $k$-connected graph of order $n$ attains
$n=\min \{n,3k+2\}.$ If $n>3k+2,$ the graph $kK_1\vee ((k+1)K_2+(n-3k-2)K_1)$ is a $k$-connected graph of order $n$ and detour order $3k+2$ which has
detour covering number $3k+2=\min \{n,3k+2\}.$ \hfill$\Box$

{\bf Definition 5.} Given a graph $H,$ a graph $G$ is said to be {\it induced-$H$ free} if $G$ contains no induced subgraph that is isomorphic to $H.$ Let $\mathscr{S}$ be a set of graphs. A graph $G$ is said to be {\it induced-$\mathscr{S}$ free} if $G$ contains no induced subgraph that is isomorphic to any graph in $\mathscr{S}.$

Next we give subgraph-forbidden conditions that ensure a graph to be detour covered or cummerbund covered. We will need the following lemma due to Golan and Shan [9].

{\bf Lemma 9 [9,} Lemma 3{\bf ].} {\it Every detour in an induced-$2K_2$ free graph is dominating.}

{\bf Theorem 10.} {\it Every connected induced-$\{P_4, 2K_2\}$ free graph is detour covered.}

{\bf Proof.} Let $G$ be a connected induced-$\{P_4,2K_2\}$ free graph and let $P=v_1v_2\dots v_t$ be a detour in $G.$ If $G$ is traceable, the conclusion holds trivially.
Next we suppose that $G$ is non-traceable. By Lemma 9, $P$ is dominating. Thus every vertex outside $P$ has a neighbor on $P,$ since $G$ is connected.
If $t\le 3,$ then $G$ is a star, which is obviously detour covered. Next we assume that $t\ge 4.$

To prove that $G$ is detour covered, it suffices to show that every vertex $v\in V(G)\setminus V(P)$ lies in a detour. Let $v_i\in N_P(v)$ for some $i.$
Note that $i\not\in \{1, t\}$ since $P$ is a detour. Since $G$ is induced-$2K_2$ free, there exists an edge $v_pv_q$ in $[\{v_1,v_2\},\,\{v_{t-1},v_t\}].$ We first observe that
  $\{p, q\}\neq \{1, t\},$ since otherwise $vv_iv_{i-1}\dots v_1v_tv_{t-1}\dots v_{i+1}$ would be a path longer than $P,$ a contradiction. There are three possibilities.

Case 1. $\{p, q\}=\{1, t-1\}.$

If $i<t-1,$ then $vv_iv_{i-1}\dots v_1v_{t-1}v_{t-2}\dots v_{i+1}$ is a detour and so $v$ lies in a detour. If $i=t-1,$ $v$ lies in the detour $vv_{t-1}v_{t-2}\dots v_1.$

Case 2. $\{p, q\}=\{2, t\}.$

This case is similar to Case 1.

Case 3. $\{p, q\}=\{2, t-1\}.$

Since $G$ is induced-$P_4$ free, the path  $v_1v_2v_{t-1}v_t$ contains a chord. Then either $v_1v_{t-1}\in E(G)$ or $v_2v_t\in E(G).$ Arguing as above, we deduce that
$v$ lies in a detour. \hfill$\Box$

{\bf Lemma 11.} {\it Every cummerbund in a $2$-connected induced-$2K_2$ free graph is dominating.}

{\bf Proof.} Let $G$ be a $2$-connected induced-$2K_2$ free graph, and let $C$ be a cummerbund in $G.$ We will show that $C$ is dominating.
To the contrary, suppose that $G-V(C)$ contains an edge $uv.$ Since $G$ is $2$-connected, there exist two vertex-disjoint $(\{u, v\}, V(C))$-paths $P$ and $Q$ where $P$ is a $(u,x)$-path and $Q$ is a $(v, y)$-path [13, p.174]. We choose $P$ and $Q$ such that  $d_C(x, y)$ is the least possible. Since $C$ is a cummerbund, we have $d_C(x, y)\ge 3.$ Let $R$
be the $(x, y)$-path on $C$ with $||R||=d_C(x, y).$ Let $wz$ be an edge of $R$ where $w$ and $z$ are internal vertices of $R.$ Since $G$ is induced-$2K_2$ free, there exists an edge in $[\{u, v\}, \{w, z\}],$ which yields two vertex-disjoint $(\{u, v\}, V(C))$-paths whose endpoints on $C$ have shorter distance on $C$ than $P, Q,$ a contradiction.
\hfill$\Box$

{\bf Theorem 12.} {\it Every $2$-connected induced-$\{P_4, 2K_2\}$ free graph is cummerbund covered.}

{\bf Proof.} Let $G$ be a $2$-connected induced-$\{P_4,2K_2\}$ free graph. If $G$ is hamiltonian, the conclusion holds trivially. Next suppose that $G$ is nonhamiltonian.
Let $C=v_1v_2\dots v_tv_1$ be a cummerbund in $G.$ To prove that $G$ is  cummerbund covered, it suffices to show that every vertex $v\in V(G)\setminus V(C)$ lies in a cummerbund.

Since $G$ is $2$-connected, by Lemma 11, $C$ is a dominating cycle of $G,$ and so $d_C(v)=d_G(v)\ge 2.$ We make the convention that $v_{t+1}=v_1.$
Let $v_i,v_j$ be two neighbors of $v$ in $C$, where $j>i$. If $d_C(v_i,v_j)=1$, then there exists a cycle longer than $C$, a contradiction.
If $d_C(v_i,v_j)=2$, then $v$ lies in a cummerbund. Next, we assume that $d_C(v_i,v_j)\ge 3.$

Since $G$ is $2K_2$-free, there exists an edge in $[\{v_{i+1},v_{i+2}\},\, \{v_{j+1},v_{j+2}\}].$ First observe that  $v_{i+1}$ and $v_{j+1}$ are nonadjacent, since otherwise
$vv_iv_{i-1}\dots v_{j+1}v_{i+1}v_{i+2}\dots v_jv$ would be a cycle longer than $C,$ a contradiction. If $v_{i+1}v_{j+2}\in E(G)$, then
$vv_iv_{i-1}\dots v_{j+2}v_{i+1}v_{i+2}\dots v_jv$ is a cummerbund containing $v.$ If $v_{i+2}v_{j+1}\in E(G)$, then $vv_iv_{i-1}\dots v_{j+1}v_{i+2}v_{i+3}\dots v_jv$ is a cummerbund containing $v.$ If $v_{i+2}v_{j+2}\in E(G)$, then since $G$ is induced-$P_4$ free, the path $v_{i+1}v_{i+2}v_{j+2}v_{j+1}$ contains a chord, which implies that $v_{i+1}v_{j+1}\in E(G)$, $v_{i+1}v_{j+2}\in E(G)$ or $v_{i+2}v_{j+1}\in E(G).$ Arguing as above, we deduce that $v$ lies in a cummerbund.
\hfill$\Box$

{\bf Remark 4.} Concerning the conditions in Theorems 10 and 12, forbidding one single induced graph in $\{P_4, 2K_2\}$ does not always ensure a graph to be detour covered or cummerbund covered. Let $a, b$ be integers with $a\ge 1$ and $b\ge 3.$ Then the $2$-connected graph $K_2\vee (aK_1+bK_2)$ is induced-$\{P_4, C_4\}$ free, but it is neither detour covered nor cummerbund covered. On the other hand, the $2$-connected graph $Q$ of order $12$ in Figure 2 is induced-$\{2K_2, C_4\}$ free, but it is neither detour covered nor cummerbund covered. The detour order and circumference of $Q$ are $11$ and $9,$ respectively, and both the detour covering number and cummerbund covering number of $Q$ are $11.$
The vertex $w$ of $Q$ does not lie in any detour or cummerbund.
\begin{figure}[h]
\centering
\includegraphics[width=0.4\textwidth]{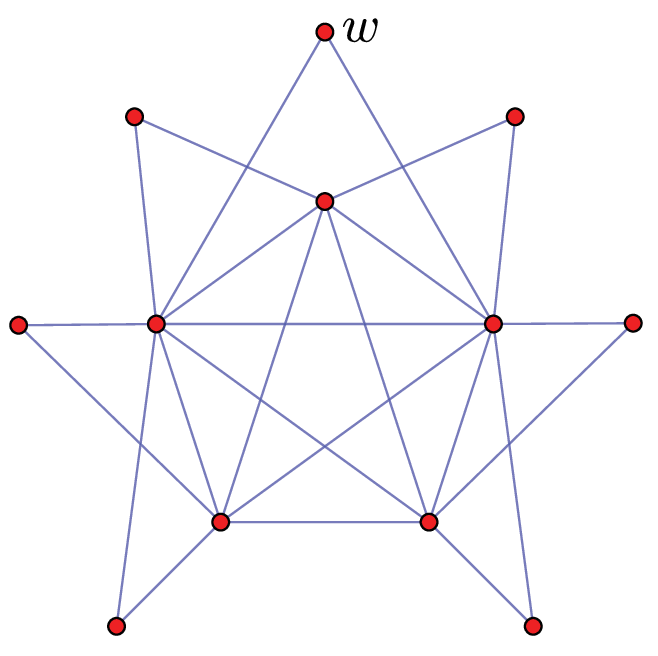}
\caption{The graph $Q$}
\end{figure}

{\bf Definition 6.} A finite simple graph $G$ is called a {\it threshold graph} if there exists a nonnegative real‐valued function $f$ defined on the vertex set of $G$, $f:V(G)\rightarrow \mathbb{R}$ and a nonnegative real number $t$ such that for any two distinct vertices $u$ and $v$, $u$ and $v$ are adjacent if and only if $f(u)+f(v)>t.$

It is known [12, p.10] that a graph is a threshold graph if and only if it is induced-$\{P_4, C_4, 2K_2\}$ free.

The following corollary follows from Theorems 10 and 12 immediately.

{\bf Corollary 13.} {\it Every connected threshold graph is detour covered and every $2$-connected threshold graph is cummerbund covered.}

Now we study bipartite graphs. Consider the seven bipartite graphs $H_i$ ($1\le i\le 7$) of order $9$ in Figure 3.
\begin{figure}[h]
\centering
\includegraphics[width=0.98\textwidth]{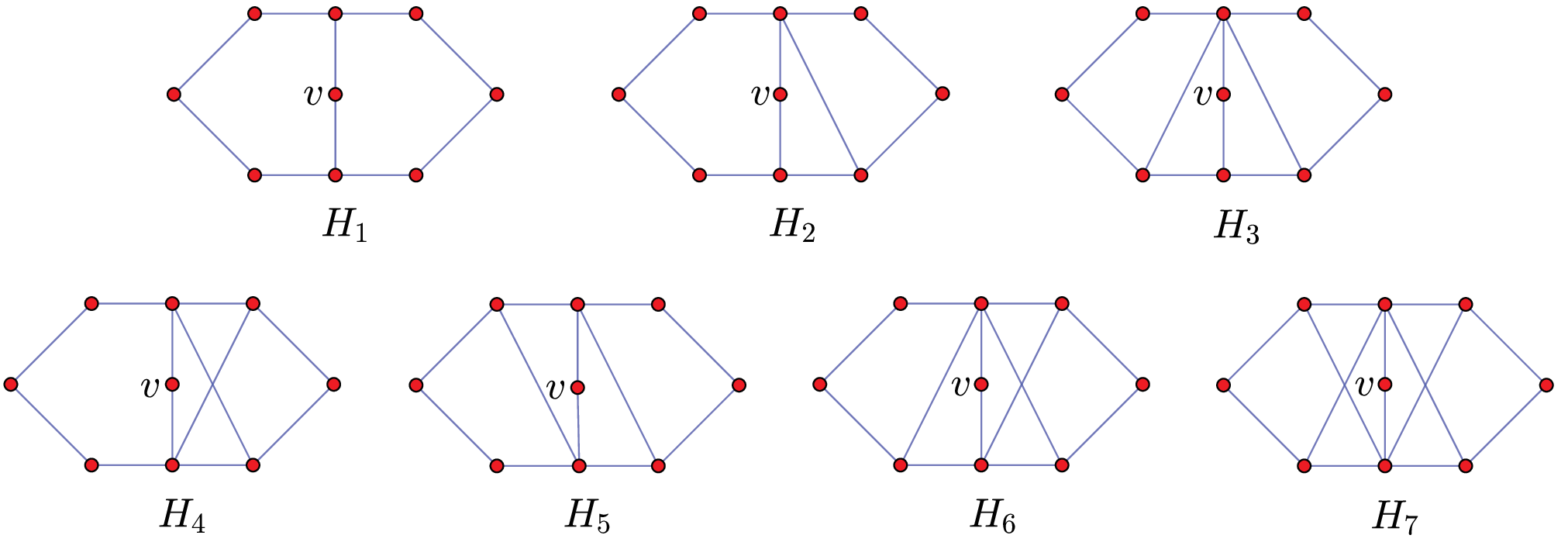}
\caption{The graphs $H_i$ $(1\le i\le 7)$}
\end{figure}
For every integer $n\ge 9$ we construct seven
bipartite graphs $G_{i,n}$ ($1\le i\le 7$) of order $n.$ First $G_{i,9}=H_i$ for $1\le i\le 7.$ For $n\ge 10$ and every $i$ with $1\le i\le 7,$ $G_{i,n}$ is obtained
from $H_i$ by duplicating the vertex $v$ $n-9$ times. Denote
$$
\mathscr{F}_n=\{G_{i,n}|\, 1\le i\le 7\}. \eqno (6)
$$

Recall that ${\rm cc}(G)$ denotes the cummerbund covering number of a graph $G.$ Given a vertex $v$ and a set $S$ of vertices, a {\it $(v, S)$-fan} is a set of
$(v, S)$-paths such that any two of them share only the vertex $v.$

{\bf Theorem 14.} {\it Let $G$ be a $2$-connected bipartite graph of order $n\ge 9$. Then ${\rm cc}(G)\ge 8.$ Equality holds if and only if $G\in \mathscr{F}_n,$ which
is defined in (6).}

{\bf Proof.} We first prove that if $c(G)\le 6,$ then $G$ is cummerbund covered, and hence ${\rm cc}(G)=n\ge 9.$
Suppose that $c(G)=6.$ Let $C$ be a cummerbund. Choose any $z\in V(G)\setminus V(C).$ Since $G$ is $2$-connected, there exists a $(z, V(C))$-fan of size two $Q_1, Q_2$
([7] or [13, p.170]), where $Q_1$ is a $(z, x)$-path and $Q_2$ is a $(z, y)$-path. Since $C$ is a cummerbund, $2\le ||Q_1||+||Q_2||\le 3$ and $2\le d_C(x, y)\le 3.$
If $d_C(x, y)=2,$ both $Q_1$ and $Q_2$ are in fact edges, and considering the subgraph $C\cup Q_1\cup Q_2$ we see that $z$ lies in a $6$-cycle. If $d_C(x, y)=3,$ we have $||Q_1||+||Q_2||=3,$ since otherwise $||Q_1||+||Q_2||=2$ and $G$ would contain a $5$-cycle, contradicting our assumption that $G$ is bipartite. Thus, in this case $z$ also lies in a $6$-cycle.

The case when $c(G)=4$ can be treated similarly as above.

Trivially ${\rm cc}(G)\ge c(G).$  If $c(G)\ge 8,$ then ${\rm cc}(G)\ge 8.$ This proves the first conclusion ${\rm cc}(G)\ge 8$ of Theorem 14.

Next we assume that ${\rm cc}(G)=8.$ Then $c(G)=8.$ Let $C=v_1v_2v_3v_4v_5v_6v_7v_8v_1$ be a cummerbund. Now every vertex in $G-V(C)$ does not lie in any cummerbund.

{\bf Claim 1.} $C$ is a dominating cycle.

To the contrary, suppose that $G-V(C)$ contains an edge $uv.$ Since $G$ is $2$-connected, there exist two vertex disjoint  $(\{u,v\},V(C))$-paths $Q_1, Q_2$
where $Q_1$ is a $(u, x)$-path and $Q_2$ is a $(v, y)$-path. Since $C$ is a cummerbund of length $8,$ and $u$ does not lie in any cummerbund, we must have
$||Q_1||=||Q_2||=1$ and $d_C(x, y)=4.$ But then the subgraph $C\cup Q_1\cup Q_2$ contains a $7$-cycle, contradicting our assumption that $G$ is bipartite.

{\bf Claim 2.} Each vertex in $G-V(C)$ has exactly two neighbors on $C$, and the distance of these two neighbors on $C$ is $4.$

By Claim 1, $G-V(C)$ is an empty graph. Let $v\in V(G)\setminus V(C).$ Since $G$ is $2$-connected, $v$ has at least two neighbors on $C.$
Let $x, y$ be any two neighbors of $v$ on $C.$ Since $G$ contains no odd cycles, $d_C(x, y)\neq 1, 3.$ Thus $d_C(x, y)\in \{2, 4\}.$ If $d_C(x, y)=2,$ $v$ would lie in
a cummerbund, a contradiction. This shows that $d_C(x, y)=4,$ which implies that $v$ has exactly two neighbors on $C$ since $C$ is an $8$-cycle. Claim 2 is proved.

A chord $xy$ of the cycle $C$ is called a {\it $k$-chord} if $d_C(x,y)=k.$ Since $G$ is bipartite and $C$ is an $8$-cycle, every chord of $C$ is a $3$-chord.
We distinguish five cases according to the number of chords of $C.$

Case 1. $C$ is chordless.

Recall that $|V(G)|\ge 9.$ Let $v\in V(G)\setminus V(C).$
By Claim 2, without loss of generality assume that $vv_1,vv_5\in E(G).$ We will show that $N_G(u)=\{v_1,v_5\}$ for each vertex $u\in V(G)\setminus V(C),$
and consequently $G=G_{1,n}$ which is defined in (6). Assume that $N_G(u)=\{u_1,u_2\}.$ By Claim 2, we have $d_C(u_1,u_2)=4.$

If $\{u_1,u_2\}=\{v_2,v_6\}$, then $v_1vv_5v_4v_3v_2uv_6v_7v_8v_1$
is a cycle of length $10$, a contradiction. If $\{u_1,u_2\}=\{v_4,v_8\},$ then $v_1vv_5v_6v_7v_8uv_4v_3v_2v_1$
is a cycle of length $10$, a contradiction. If $\{u_1,u_2\}=\{v_3,v_7\}$, then $u$ lies in the cummerbund $v_1vv_5v_4v_3uv_7v_8v_1,$
a contradiction. Therefore, $\{u_1,u_2\}=\{v_1,v_5\}.$

Case 2. $C$ has exactly one chord.

Without loss of generality, we assume that $v_1v_4\in E(G).$
Let $u$ be any vertex in $V(G)\setminus V(C)$ and let $N_G(u)=\{u_1,u_2\}.$
By Claim 2, we have $d_C(u_1,u_2)=4.$ If $\{u_1,u_2\}=\{v_2,v_6\}$, then $u$ lies in the cummerbund $v_1v_4v_3v_2uv_6v_7v_8v_1,$
a contradiction.
If $\{u_1,u_2\}=\{v_3,v_7\}$, then $u$ lies in the cummerbund $v_1v_2v_3uv_7v_6v_5v_4v_1,$
a contradiction.
Thus, $\{u_1,u_2\}=\{v_1,v_5\}$ or $\{u_1,u_2\}=\{v_4,v_8\}.$ If there exist two vertices $w_1,w_2\in V(G)\setminus V(C)$ such that
$N_G(w_1)=\{v_1, v_5\}$ and $N_G(w_2)=\{v_4, v_8\},$ then $v_1v_2v_3v_4w_2v_8v_7v_6v_5w_1v_1$ is a cycle of length $10$, a contradiction.
Hence either $N_G(w)=\{v_1,v_5\}$ for each vertex $w\in V(G)\setminus V(C)$ or  $N_G(w)=\{v_4,v_8\}$ for each vertex $w\in V(G)\setminus V(C).$
Thus $G=G_{2,n}.$

Case 3. $C$ has exactly two chords.

Essentially there are four possibilities for the two chords to lie in $C.$
Without loss of generality, we assume that $v_1v_4,v_1v_6\in E(G),$ or $v_1v_4,v_2v_5\in E(G),$ or $v_1v_4,v_5v_8\in E(G),$ or $v_1v_4,v_2v_7\in E(G).$

Suppose that $v_1v_4,v_1v_6\in E(G).$ Let $u\in V(G)\setminus V(C)$. By the proof of Case 2, we have $N_G(u)=\{v_1,v_5\}$ and so $G=G_{3,n}$.

Suppose that $v_1v_4, v_2v_5\in E(G).$ Let $u\in V(G)\setminus V(C)$. By the proof of Case 2, we have $N_G(u)=\{v_1,v_5\}$ and so $G=G_{4,n}$.

Suppose that $v_1v_4,v_5v_8\in E(G).$ Let $u\in V(G)\setminus V(C)$. By the proof of Case 2, we have $N_G(u)=\{v_1,v_5\}$ or $N_G(u)=\{v_4,v_8\}$.
If there exist two vertices $w_1,w_2\in V(G)\setminus V(C)$ such that $N_G(w_1)=\{v_1, v_5\}$ and $N_G(w_2)=\{v_4, v_8\},$ then $v_1v_2v_3v_4w_2v_8v_7v_6v_5w_1v_1$
is a cycle of length $10$, a contradiction. Hence either $N_G(w)=\{v_1,v_5\}$ for each vertex $w\in V(G)\setminus V(C)$ or  $N_G(w)=\{v_4,v_8\}$ for each vertex
$w\in V(G)\setminus V(C).$ Thus $G=G_{5,n}.$

Suppose that $v_1v_4,v_2v_7\in E(G).$ Let $u\in V(G)\setminus V(C).$ According to the proof of Case 2, by the existence of the chord $v_1v_4$, we have $N_G(u)=\{v_1,v_5\}$ or $N_G(u)=\{v_4,v_8\}.$ However, by the existence of the chord $v_2v_7$, we have $N_G(u)=\{v_2,v_6\}$ or $N_G(u)=\{v_3,v_7\}$, a contradiction.

Case 4. $C$ has exactly three chords.

According to the proof of Case 3, we may assume that $v_1v_4,v_1v_6,v_2v_5\in E(G)$. Furthermore, we have $N_G(u)=\{v_1,v_5\}$ for each vertex $u\in V(G)\setminus V(C).$
Thus $G=G_{6,n}.$

Case 5. $C$ has exactly four chords.

According to the proof of Case 3, we may assume that $v_1v_4,v_1v_6,v_2v_5,v_5v_8\in E(G)$. Furthermore, we have $N_G(u)=\{v_1,v_5\}$ for each vertex $u\in V(G)\setminus V(C).$
Thus $G=G_{7,n}.$

Case 6. $C$ has at least five chords.

According to the proof of Case 3, this case cannot occur.\hfill$\Box$

The {\it theta graph} $\theta(l_1, l_2, \dots, l_k)$ is the graph obtained by connecting two vertices by pair-wise internally disjoint paths of lengths $l_1, l_2, \dots, l_k.$
We write $\theta(a_1^{m_1}, a_2^{m_2}, \dots, a_s^{m_s})$ for
$\theta(\operatorname*{\underbrace{a_1\cdots a_1}}\limits_{m_1}, \operatorname*{\underbrace{a_2\cdots a_2}}\limits_{m_2},\ldots, \operatorname*{\underbrace{a_s\cdots a_s}}\limits_{m_s}).$ A theta graph is called {\it uniform} if it is of the form $\theta(a^m).$ Next we consider girth conditions for the cummerbund covering number.
The following lemma must be known, but we cannot locate a reference.

{\bf Lemma 15.} {\it A $2$-connected graph $G$ satisfies $g(G)=c(G)$ if and only if $G$ is a cycle or a uniform theta graph.}

{\bf Proof.} The ``if part" is trivial. We prove the ``only if part." Let $G$ be a $2$-connected graph with $g(G)=c(G).$
The condition $g(G)=c(G)$ means that all cycles in $G$ have the same length. Denote $k=g(G).$ Let $C$ be a cycle.

If $k$ is odd, we assert that $G=C.$ Otherwise there is a vertex $x$ outside $C.$ Since $G$ is $2$-connected, there is a $(x, V(C))$-fan consisting of two paths $P_1$
and $P_2$ ([7] or [13, p.170]). Then the subgraph $C\cup P_1\cup P_2$ contains a theta graph, which contains an even cycle, a contradiction.

Now suppose $k$ is even. Let $k=2p.$ Suppose $G$ is not a cycle. Choose a vertex $z\in V(G)\setminus V(C).$ Then since $G$ is $2$-connected, there is a $(z, V(C))$-fan consisting of two paths $Q_1$ and $Q_2,$ where $Q_1$ is a $(z,x)$-path and $Q_2$ is a $(z,y)$-path. Denote $Q=Q_1\cup Q_2.$ Since all cycles in $G$ have the same length $2p,$ the theta graph $T$ consisting of the three paths $C(x, y), C(y, x)$ and $Q$ must be a uniform theta graph; i.e., $d_C(x,y)=p=||Q||.$ Clearly $G[V(C)\cup V(Q)]=T,$ since otherwise $G$ would contain a
shorter cycle. If $V(G)\setminus V(T)\neq\emptyset,$ let $w\in V(G)\setminus V(T).$ Arguing as above, we deduce that there exists a $(u, v)$-path $R$ of length $p$ containing $w$  that is internal disjoint from $C$ where $u,v\in V(C)$ and $d_C(u,v)=p$ such that  $G[V(C)\cup V(R)]$ is a uniform theta graph.

Denote $\Omega=(V(Q)\cap V(R))\setminus\{x,y,u,v\}.$ We assert that $Q$ and $R$ are internally disjoint; i.e., $\Omega=\emptyset.$ To the contrary, suppose $\Omega\neq\emptyset.$
It is easy to verify that if $|\Omega|=1,$ then $C\cup Q\cup R$ contains a cycle of length less than $2p,$ a contradiction, and if $|\Omega|\ge 2,$ then $Q\cup R$ contains a cycle of length less than $2p,$ a contradiction again.

Next we assert that $\{u, v\}=\{x, y\};$ i.e., $R$ and $Q$ have the same endpoints. Otherwise  $C\cup Q\cup R$ would contain a cycle of length greater than $2p,$ a contradiction.
So far, we have shown that $C\cup Q\cup R$ is a uniform theta graph. Continuing the above process we conclude that $G$ is a uniform theta graph.\hfill$\Box$

{\bf Theorem 16.} {\it The minimum cummerbund covering number of a $2$-connected graph of order at least $8$ with girth at least $4$ is $6.$
The minimum cummerbund covering number of a $2$-connected graph of order at least $10$ with girth at least $5$ is $8.$
The minimum cummerbund covering number of a $2$-connected graph of an even order at least $12$ with girth at least $6$ is $8,$ and the minimum cummerbund covering number of a $2$-connected graph of an odd order at least $13$ with girth at least $6$ is $9.$}

{\bf Proof.} (1) Let $G$ be a $2$-connected graph of order at least $6$ with $g(G)\ge 4$. We first show ${\rm cc}(G)\ge 6$.
To the contrary, we assume ${\rm cc}(G)\le 5.$ Clearly, $c(G)\le 5.$

Let $C$ be a cummerbund. Choose any vertex $v$ outside $C.$ Since $G$ is $2$-connected, there is a $(v, V(C))$-fan consisting of
two paths $P_1$ and $P_2$ ([7] or [13, p.170]). Since  $c(G)\le 5$ and $C$ is a cummerbund,  $||P_1||=||P_2||=1,$ which implies that
$v$ lies in a cummerbund. Since $v$ was arbitrarily chosen, $G$ is cummerbund covered, a contradiction.

Conversely, for every integer $n\ge 8,$ the theta graph $\theta(3^2,2^{n-6})$ is a $2$-connected graph of order $n$ with girth $4$ whose cummerbund covering number
is $6.$

(2) Let $G$ be a $2$-connected graph of order at least $10$ with $g(G)\ge 5.$ We first prove ${\rm cc}(G)\ge 8$.
To the contrary, assume that ${\rm cc}(G)\le 7.$ By Lemma 15, if $c(G)=5$ then $G=C_5,$ contradicting our assumption that $G$ has order at least $10.$
Hence  $6\le c(G)\le 7$.

Case 1. $c(G)=7.$

Let $C=v_1v_2v_3v_4v_5v_6v_7v_1$ be a cummerbund. Let $z\in V(G)\setminus V(C)$. Since $G$ is $2$-connected, there is a $(z, V(C))$-fan consisting of
two paths $P_1$ and $P_2.$ Let $\{u_i\}=V(P_i)\cap V(C)$ for $i=1,2$. Since $C$ is a cummerbund, we have $d_C(u_1,u_2)\ge 2.$
Suppose that $d_C(u_1,u_2)= 2.$ Without loss of generality, assume that $u_1=v_1$ and $u_2=v_3$. Since $g(G)\ge 5$, we have $||P_1||\ge 2$ or $||P_2||\ge 2$ and so $v_1P_1zP_2v_3v_4v_5v_6v_7v_1$ is a cycle of length at least $8$, a contradiction.

Suppose that $d_C(u_1,u_2)= 3.$ Without loss of generality, assume that $u_1=v_1$ and $u_2=v_4$. If $||P_1||\ge 2$ or $||P_2||\ge 2,$ then $||P_1||+||P_2||\ge 3$.
This implies that $z$ lies in a cummerbund and so ${\rm cc}(G)\ge 8$, a contradiction. Hence  $||P_1||=||P_2||=1;$ i.e., $zv_1\in E(G)$ and $zv_4\in E(G)$.
Since $|V(G)|\ge 10,$ there exists a vertex $y\in V(G)\setminus (V(C)\cup \{z\})$. Since $G$ is $2$-connected, there exists a $(y,V(C))$-fan consisting of two paths
$Q_1, Q_2$. Let $\{w_i\}=V(Q_i)\cap V(C)$ for $i=1, 2.$ Similarly as above, we only need to consider the case when  $d_C(w_1,w_2)=3$ and $||Q_1||=||Q_2||=1.$
Since $g(G)\ge 5,$ we have $\{w_1,w_2\}\neq \{v_1,v_4\}$. If $\{w_1,w_2\}=\{v_2,v_5\}$, then $v_1zv_4v_3v_2yv_5v_6v_7v_1$ is a cycle of length $9$, a contradiction. If $\{w_1,w_2\}=\{v_3,v_6\}$, then $v_1zv_4v_3yv_6v_7v_1$ is a cummerbund, and so $z$ lies in a cummerbund. Hence, ${\rm cc}(G)\ge 8$, a contradiction. If $\{w_1,w_2\}=\{v_6,v_2\}$ or $\{w_1,w_2\}=\{v_7,v_3\}$, then there exists a cycle of length at least $8$, a contradiction. It follows that $\{w_1,w_2\}=\{v_1,v_5\}$ or $\{w_1,w_2\}=\{v_4,v_7\}.$

Without loss of generality, we assume that $uv_1, uv_5\in E(G)$. Since $|V(G)|\ge 10,$ there exists a vertex $x\in V(G)\setminus(V(C)\cup \{z,y\})$.
Similarly as above, we have $xv_1, xv_5\in E(G)$ or $xv_4, xv_7\in E(G).$ If $xv_1, xv_5\in E(G)$, then $v_1xv_5yv_1$ is a cycle of length $4$, a contradiction.
If $xv_4,xv_7\in E(G)$, then $v_1yv_5v_6v_7xv_4v_3v_2v_1$ is a cycle of length $9$, a contradiction.

Case 2. $c(G)=6.$

Let $C=v_1v_2v_3v_4v_5v_6v_1$ be a cummerbund. Let $v\in V(G)\setminus V(C)$. Since $G$ is $2$-connected, there exists a $(v,V(C))$-fan consisting of two paths
$P_1, P_2$. Let $\{u_i\}=V(P_i)\cap V(C)$ for $i=1,2$. Clearly, $d_C(u_1,u_2)\ge 2$.

Suppose that $d_C(u_1,u_2)=2$. Without loss of generality, we assume that $u_1=v_1, u_2=v_3.$ Since $g(G)\ge 5$, we have $||P_1||\ge 2$ or $||P_2||\ge 2$ which implies that $v_1P_1vP_2v_3v_4v_5v_6v_1$ is a cycle of length at least $7$, a contradiction.

Suppose that $d_C(u_1,u_2)=3$. Without loss of generality, we assume that $u_1=v_1, u_2=v_4.$ If $||P_1||\ge 2$ or $||P_2||\ge 2$, then $P_1\cup P_2$ is a path
contained in a cummerbund, and so ${\rm cc}(G)\ge 8,$ a contradiction. Thus $||P_1||=||P_2||=1;$ i.e., $vv_1,vv_4\in E(G).$
Let $w$ be a vertex in $V(G)\setminus (V(C)\cup \{v\})$. Let $Q_1,Q_2$ be a $(w,V(C))$-fan and let $\{w_i\}=V(Q_i)\cap V(C)$ for $i=1,2.$
As above, we have $d_C(w_1,w_2)=3$ and $||Q_1||=||Q_2||=1$.

If $\{w_1,w_2\}=\{v_1,v_4\}$, then $v_1wv_4vv_1$ is a cycle of length $4$, a contradiction. If $\{w_1,w_2\}=\{v_2,v_5\}$, then $v_1vv_4v_3v_2wv_5v_6v_1$ is a cycle of length $8$, a contradiction. If $\{w_1,w_2\}=\{v_3,v_6\}$, then $v_1vv_4v_5v_6wv_3v_2v_1$ is a cycle of length $8$, a contradiction.

Conversely, if $n\ge 11$ is an odd integer, let $H=\theta(4^2,3^{(n-9)/2},2).$  If $n\ge 10$ is an even integer, let $xpqry$ be a path in the theta graph $T=\theta(4^2, 3^{(n-10)/2}, 2)$ where $x, y$ are the two vertices with degree greater than $2.$ Now construct a graph $H$ from $T$ by first adding the edge $xr$ and then subdividing it. Then $H$ is a $2$-connected graph of order $n$ with girth $5$ whose cummerbund covering number is $8.$

(3) Let $G$ be a $2$-connected graph of order $n\ge 12$ with $g(G)\ge 6.$ By part (2) we have ${\rm cc}(G)\ge 8$.

We first prove that if $n\ge 13$ is odd then ${\rm cc}(G)\ge 9.$ To the contrary, assume that ${\rm cc}(G)\le 8.$ By Lemma 15, if $c(G)=6$ then $G$ is a uniform theta graph and so $G$ is cummerbund covered, a contradiction. Hence $7\le c(G)\le 8.$

Case 1. $c(G)=8$.

Let $C=v_1v_2v_3v_4v_5v_6v_7v_8v_1$ be a cummerbund. Since $G$ is $2$-connected, for any vertex $u\in V(G)\setminus V(C),$ there exists a $(u, V(C))$-fan consisting of two paths $P_1,P_2.$ We choose a vertex $u$ and two paths $P_1,P_2$ such that $||P_1||+||P_2||$ is maximized.  Let $\{u_i\}=V(P_i)\cap V(C)$ for $i=1,2$.

Subcase 1.1. $||P_1||=||P_2||=1$.

 Since $g(G)\ge 6,$ we have $d_C(u_1,u_2)=4.$ Without loss of generality, we assume that $u_1=v_1,u_2=v_5.$
Since $|V(G)|\ge 13$, there exists a vertex $w\in V(G)\setminus (V(C)\cup \{u\})$. Since $G$ is $2$-connected, there exists a $(w,V(C))$-fan consisting of two paths $Q_1,Q_2.$ Let $\{w_i\}= V(Q_i)\cap V(C)$ for $i=1,2.$ By the choice of $P_1,P_2$, we have $||Q_1||=||Q_2||=1,$ and so $d_C(w_1,w_2)=4.$
If $\{w_1,w_2\}=\{v_1,v_5\}$, then $v_1uv_5wv_1$ is a cycle of length $4$, a contradiction.
If $\{w_1,w_2\}=\{v_2,v_6\}$, then $v_1uv_5v_4v_3v_2wv_6v_7v_8v_1$ is a cycle of length $10$, a contradiction. If $\{w_1,w_2\}=\{v_4,v_8\},$ then $v_1uv_5v_6v_7v_8wv_4v_3v_2v_1$ is a cycle of length $10$, a contradiction. If $\{w_1,w_2\}=\{v_3,v_7\}$, then $u$ lies in a cummerbund $v_1uv_5v_4v_3wv_7v_8v_1$ and so ${\rm cc}(G)\ge 9,$ a contradiction.

Subcase 1.2. $\{||P_1||,||P_2||\}=\{1,2\}.$

In this case, we have $d_C(u_1,u_2)=4$, since otherwise $u$ lies in a cummerbund and so ${\rm cc}(G)\ge 9$, a contradiction.
Without loss of generality, we assume that $P_1=uvv_1$ and $P_2=uv_5.$ Since $|V(G)|\ge 13,$ there exists a vertex $w\in V(G)\setminus (V(C)\cup \{u,v\}).$ Since $G$ is $2$-connected, there exists a $(w,V(C))$-fan consisting of two paths $Q_1,Q_2.$ Let $\{w_i\}= V(Q_i)\cap V(C)$ for $i=1,2$.

We assert that $\{u,v\}\cap (V(Q_1)\cup V(Q_2))=\emptyset.$ To the contrary, we assume that $v\in V(Q_1)\cup V(Q_2).$ By the choice of $P_1,P_2,$ we have $vw,v_5w\in E(G).$ However, $v_5uvwv_5$ would be a $4$-cycle, a contradiction.

Next we assert that $d_C(w_1,w_2)=4.$ By the choice of $P_1,P_2,$ we have either $||Q_1||=||Q_2||=1$ or $\{||Q_1||, ||Q_2||\}=\{1,2\}.$ If $||Q_1||=||Q_2||=1,$ then we have $d_C(w_1,w_2)=4,$ since otherwise $G$ would contain a cycle of length at most $5$, a contradiction. If $\{||Q_1||, ||Q_2||\}=\{1,2\},$ we also have $d_C(w_1,w_2)=4,$ since otherwise
$w$ would lie in a cummerbund and so ${\rm cc}(G)\ge 9$, a contradiction.

If $\{w_1,w_2\}=\{v_2,v_6\},$ then $v_1vuv_5v_4v_3v_2Q_1wQ_2v_6v_7v_8v_1$ is a cycle of length at least $11$ and so $c(G)\ge 11,$ a contradiction. If $\{w_1,w_2\}=\{v_4,v_8\},$ then $v_1vuv_5v_6v_7v_8Q_1wQ_2v_4v_3v_2v_1$ is a cycle of length at least $11,$ a contradiction. If  $\{w_1,w_2\}=\{v_3,v_7\},$ then $G$ contains a cycle
$v_1vuv_5v_4v_3Q_1wQ_2v_7v_8v_1$ of length at least $9,$  a contradiction. Thus, we have $\{w_1,w_2\}=\{v_1,v_5\}.$ Since $g(G)\ge 6,$ we have $\{||Q_1||, ||Q_2||\}=\{1,2\},$ otherwise $wv_1vuv_5w$ would be a $5$-cycle, a contradiction. Hence, $T=C\cup P_1\cup P_2\cup Q_1\cup Q_2$ is a theta graph $\theta(4^2,3^2).$ Moreover, $T$ is an induced subgraph of $G$, since otherwise $G$ would contain a cycle of length less than $6$, a contradiction. Continuing the above process we conclude that $G$ is a theta graph $\theta(4^2,3^{(n-8)/2}),$ contradicting our assumption that $G$ has odd order.

Subcase 1.3. $||P_1||\ge 2$ and $||P_2||\ge 2.$

Then $u_1P_1uP_2u_2Cu_1$ is a cycle of length at least $8,$ where $u_2Cu_1$ is a $(u_2,u_1)$-path on $C$ of length at least $4$. Thus $u$ lies in a cummerbund and so ${\rm cc}(G)\ge 9$, a contradiction.

Case 2. $c(G)=7.$

Let $C$ be a cummerbund. Let $u\in V(G)\setminus V(C).$ Since $G$ is $2$-connected, there  exists a $(u,V(C))$-fan consisting of two paths $P_1,P_2.$ Let $\{u_i\}=V(P_i)\cap V(C)$ for $i=1,2.$

If $||P_1||\ge 2$ and $||P_2||\ge 2,$ then $u_1P_1uP_2u_2Cu_1$ is a cycle of length at least $8$, where $u_2Cu_1$ is a $(u_2,u_1)$-path on $C$ of length at least $4$, a contradiction.
If $||P_1||=||P_2||=1,$ then there exists a cycle $u_1Cu_2uu_1$ of length at most $5$, where $u_1Cu_2$ is a $(u_1,u_2)$-path on $C$ of length at most $3$, a contradiction.

Suppose that $\{||P_1||,||P_2||\}=\{1,2\}.$ Since $g(G)\ge 6,$ we have $d_C(u_1,u_2)=3.$ Then $P_1\cup P_2$ is a path contained in a cummerbund $u_1P_1uP_2u_2Cu_1,$ where $u_2Cu_1$ is a $(u_2,u_1)$-path on $C$ of length $4$, and so ${\rm cc}(G)\ge 9,$ a contradiction.

So far, we have proved that if $n\ge 12$ is even then ${\rm cc}(G)\ge 8,$  and if $n\ge 13$ is odd then ${\rm cc}(G)\ge 9.$

Conversely, if $n\ge 12$ is an even integer, let $H=\theta(4^2,3^{(n-8)/2}).$ Then $H$ is a $2$-connected graph of order $n$ with girth $6$ whose cummerbund covering number is $8.$
If $n\ge 13$ is an odd integer, let $H=\theta(5,4,3^{(n-9)/2}).$ Then $H$ is a $2$-connected graph of order $n$ with girth $6$ whose cummerbund covering number is $9.$
\hfill$\Box$

Finally we pose several unsolved problems.

{\bf Problem 1.} Let $g\ge 7$ be an integer. Determine the minimum cummerbund covering number of a $2$-connected graph of order $n$ and girth at least $g.$

{\bf Problem 2.} Determine the minimum cummerbund covering number of a $2$-connected claw-free graph of order $n.$

{\bf Problem 3.} Determine the minimum detour covering number of a connected claw-free graph of order $n.$

\vskip 5mm
{\bf Acknowledgement.} This research  was supported by the NSFC grant 12271170 and Science and Technology Commission of Shanghai Municipality
 grant 22DZ2229014.

\end{document}